\documentclass[a4paper,12pt]{article}
  \usepackage{amssymb,amsbsy,amsmath,latexsym,theorem,eepic,enumerate,graphicx}
  \usepackage{xypic,times,overpic} 
  \xyoption{curve}
  \usepackage{geometry}
  \usepackage[ps2pdf, bookmarks=false, colorlinks=false,
            linkcolor=Green, citecolor=Red,
            urlcolor=Blue, filecolor=Green,
            pagecolor=Green,
            bookmarksnumbered=true, breaklinks=true,
            pdfstartview=FitH, hyperfigures=false,
            plainpages=false, naturalnames=true,
            pdfpagelabels, pagebackref=false]{hyperref}
  \usepackage[usenames]{color}
{\theorembodyfont{\rmfamily}\newtheorem{example}{Example}[section]}
{\theorembodyfont{\rmfamily}\newtheorem{Def}[example]{Definition}}
\newtheorem{prop}[example]{Proposition}
\newtheorem{thm}[example]{Theorem}
{\theorembodyfont{\rmfamily}\newtheorem{rem}[example]{Remark}}
\newtheorem{cor}[example]{Corollary}
\newtheorem{lem}[example]{Lemma}

  \def\FTop{\mathsf{FTop}}

\def\Set{\mathsf{Set}}

\def\A{\alpha}
\def\geq{\geqslant}
\def\leq{\leqslant}

\newcommand{\labto}[1]{\stackrel{#1}{\longrightarrow}}

\newenvironment{proof}{\noindent {\bf Proof }}{\rule{0mm}{1mm}\hfill $\Box$

\mbox{}}

\newcommand{\threeaxes}[3]{\def\objectstyle{\scriptstyle}  \objectmargin={0pt}
\xy
(0,0)*+{}="a",(0,-6)*+{\rule{0em}{1.5ex}#2}="b",(7,0)*+{\;#1}="c",
(14,-3)*+{\;#3}="d" \ar@{->} "a";"b" \ar @{->}"a";"c"  \ar
@{->}"a";"d"\endxy }

\newcommand{\directs}[2]{\def\objectstyle{\scriptstyle}  \objectmargin={0pt}
\xy
(0,4)*+{}="a",(0,-2)*+{\rule{0em}{1.5ex}#2}="b",(7,4)*+{\;#1}="c"
\ar@{->} "a";"b" \ar @{->}"a";"c" \endxy }

\newcommand{\xdirects}[2]{\def\objectstyle{\scriptstyle}  \objectmargin={0pt}
\xy
(0,0)*+{}="a",(0,-6)*+{\rule{0em}{1.5ex}#2}="b",(7,0)*+{\;#1}="c"
\ar@{->} "a";"b" \ar @{->}"a";"c" \endxy }
\newcommand{\sdirects}[2]{\def\objectstyle{\scriptstyle}  \objectmargin={0pt}
\xy
(0,2.2)*+{}="a",(0,-2.5)*+{\rule{0em}{1.5ex}#2}="b",(7,2.2)*+{\;#1}="c"
\ar@{->} "a";"b" \ar @{->}"a";"c" \endxy }

 \def\I{{\Box}\;}

\def\del{\partial}
\def\A{\alpha}
\def\B{\beta}
\def\eps{\varepsilon}
\def\epsilon{\varepsilon}

\def\half{\frac{1}{2}}

\def\cal{\mathcal}
\def\le{\leqslant}
\def\ge{\geqslant}
\def\z{\mathbf{z}}

\def\x{\mathbf{x}}

\def\d{\mathbf{d}}

\def\del{\partial}
\def\A{\alpha}
\def\B{\beta}
\def\eps{\varepsilon}
\def\epsilon{\varepsilon}

\def\half{\frac{1}{2}}

\def\cal{\mathcal}
\def\le{\leqslant}
\def\ge{\geqslant}
\def\leq{\leqslant}
\def\geq{\geqslant}
\def\subset{\subseteq}
 \def\c{\mathbin{\#}}

 \def\d{\partial}

 \def\Im{\mathop{\rm Im}\nolimits}

 \def\R{{\mathbb R}}
 \def\eps{\varepsilon}

\def\A{\alpha}

\def\ogpd{$\omega$-$\mathsf{Gpd}$}

\def\FTop{\mathsf{FTop}}

\def\cal{\mathcal}
\def\epsilon{\varepsilon}

\def\subs{\subseteq}

\def\disk{\xy*\cir<3.5pt>{}\endxy}

\def\bigcirc{{\disk}}
  
\begin{document}

  \title{A new higher homotopy groupoid: \\the fundamental  globular $\omega$-groupoid \\
  of a filtered space}
  \author{Ronald Brown\thanks{School of Computer Science, University of Wales, Dean St., Bangor, Gwynedd, LL57 1UT, UK;
email: r.brown@bangor.ac.uk}\;\thanks{The author was supported for
part of this work by a Leverhulme Emeritus Fellowship (2002-2004).
}} \maketitle
\begin{abstract}
We use the  $n$-globe with its skeletal filtration to define the
fundamental globular $\omega$--groupoid of a filtered space; the
proofs use an analogous fundamental cubical $\omega$--groupoid due
to the author and Philip Higgins. This method also relates the
construction to the fundamental crossed complex of a filtered space,
and this relation allows the proof that the crossed complex
associated to the free globular $\omega$-groupoid on one element of
dimension $n$ is the fundamental crossed complex of the $n$-globe.
\end{abstract}
{\small \tableofcontents }
\section*{Introduction}\label{sec:int}
\addcontentsline{toc}{section}{Introduction} By the $n$-globe $G^n$
we mean the subspace of Euclidean $n$-space $\R^n$ of points $x$
such that $\|x\|\leq 1$ but with the cell structure for $n \geq 1$
\begin{equation}
  G^n= e^0_\pm \cup e^1_\pm \cup \cdots \cup e^{n-1}_\pm \cup e^n.
\end{equation}This structure will be given precisely in section
\ref{sec:disksglobes}.

\begin{center}
\begin{overpic}[scale=0.7]{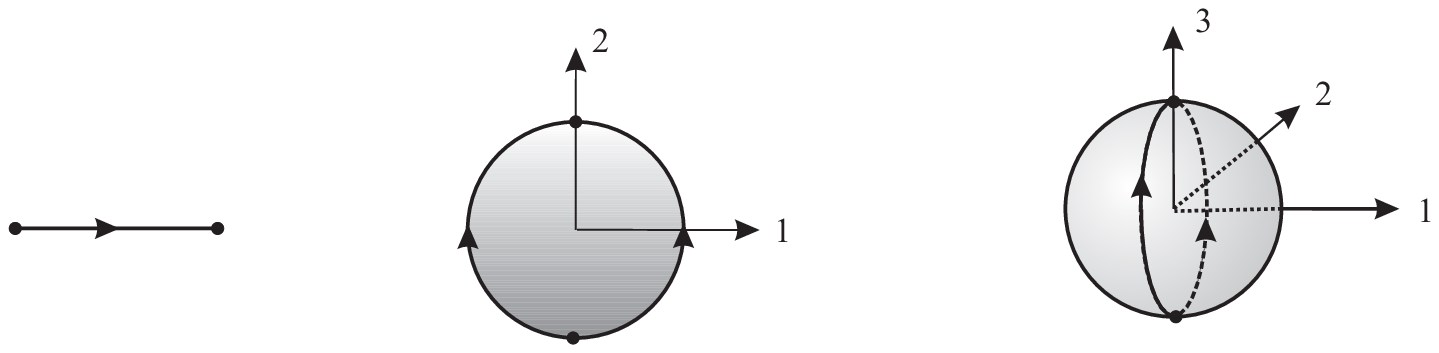}
\end{overpic}
\mbox{ }\\
\mbox{ }\\
\end{center}

A {\it filtered space } is a compactly generated space $X_\infty$
and a sequence of subspaces
\begin{equation}
X_*: X_0 \subs X_1 \subs \cdots \subs X_n \subs \cdots \subs
X_\infty .
\end{equation}
A {\it map of filtered spaces} $f: Y_* \to X_*$ is a map $f:
Y_\infty \to X_\infty$ such that $f(Y_n) \subs X_n$ for all $ n \geq
0$. This gives the category $\FTop$ of filtered spaces. A {\it
filter homotopy} $f_t: f_0 \simeq f_1$ is a continuous family of
filtered maps $f_t: Y_* \to X_*$ for $0 \leq t \leq 1$.

The $n$-globe $G^n$ has a skeletal filtration giving a filtered
space $G^n_*$. If $X_*$ is a filtered space then we have a {\it
globular singular complex} $R^\bigcirc\, X_*$ which in dimension $n$
is $\FTop(G^n_*, X_*)$. We will in appendix  \ref{sec:globular}
explain the structure of $R^\bigcirc\, X_*$ as a {\it globular set}.

We define \begin{equation}\rho^\bigcirc X_*  =(R^\bigcirc \,X_*/
\equiv),
\end{equation}
where $\equiv$ is the relation of filter homotopy {\it rel
vertices}. It will be clear that $ \rho^\bigcirc \,X_*$ inherits
from $R^\bigcirc X_* $ the structure of globular set. Our main
result is the following:
\begin{thm}[Main theorem]
There are   compositions $\circ_i, 1 \leq i \leq n$ in dimensions
$n\geq 1$ giving the globular set $ \rho^\bigcirc \, X_*$ the
structure of globular $\omega$-groupoid.
\end{thm}
We call $ \rho^\bigcirc\, X_*$ the {\it fundamental globular higher
homotopy groupoid of the filtered space $X_*$}. The proof of this
theorem goes via the notion of cubical higher homotopy groupoid of a
filtered space, established in \cite{BH81:col}. It should be useful
therefore to put these results in context.

The overall aim of work on higher homotopy groupoids may be subsumed
in the following diagram and its properties:
\begin{equation}
 \def\labelstyle{\textstyle}\vcenter{ \xymatrix@R=3pc{ \txt{\rm topological
 data}
  \ar @<0.5ex>[rr] ^\Xi   \ar [dr]_U && \txt{algebraic data}  \ar @<0.5ex>[ll] ^{\Bbb{B}} \ar [dl] ^B \\
&\txt{\rm topological spaces} }}
\end{equation}The aim is to find suitable categories of topological data,
algebraic data and functors as above, where $U$ is the forgetful
functor and $B= U \circ \Bbb B$, with the following properties:
\begin{enumerate}[(1)]
\item the functor $\Xi$ is defined homotopically and satisfies a higher
homotopy van Kampen  theorem (HHvKT)\footnote{Jim Stasheff has
suggested this term to the author, instead of the previously used
Generalised van Kampen Theorem, to make clear the higher homotopy
information contained in theorems of this type. }, in that it
preserves certain colimits;
\item $\Xi \circ \Bbb B$ is naturally equivalent  to $1$;
\item there is a natural transformation $1 \to {\Bbb B} \circ \Xi$
preserving some homotopical information.
\end{enumerate}
The purpose of (1) is to allow some calculation of $\Xi$ by gluing
simple examples, such as convex subsets, following the use of the
fundamental groupoid in \cite{B:book1}. This condition (1) at
present also rules out  some widely used algebraic data, such as for
example simplicial groups or groupoids, or differential graded
algebras, since for those cases no such functor $\,\Xi\,$ is known.
(2) shows that the algebraic data faithfully captures some of the
topological data. The imprecise (3) gives further information on the
algebraic modelling. The functor $B$ should be called a {\it
classifying space} because it often generalises the classifying
space of a group or groupoid. It has also been found useful in the
homotopy classification of maps.

Here is a table illustrating the possibilities.
\begin{center}

\begin{tabular}{|c|c|}
\hline {\it Topological data} & {\it Algebraic data}
\\ \hline \hline
space with base point & groups \\
\hline space with set of base points & groupoids \\\hline pointed
pair of spaces & crossed modules \\\hline filtered space & crossed
complexes\\\hline $n$-cube of pointed spaces& cat$^n$-groups \\
\hline $n$-cube of pointed spaces & crossed $n$-cube of groups\\
\hline
\end{tabular}
\end{center}

Strong  results in the last two cases are shown in \cite{BL87,ESt}.

In this paper we will deal only with the  case of filtered spaces,
which of course includes the first three cases. There are still
several choices of algebraic data as shown in the following diagram
of equivalent categories, which is taken from \cite{B:gpds-crossed}:
\begin{equation}    \label{diag:manycats} \def\labelstyle{\textstyle}\vcenter{\xymatrix@R-2ex@C=3pc{
   \txt{cubical \\ $ T$-complexes} \ar[r] ^-{(a)}&\txt{cubical \\$\omega$-groupoids \\with connections}
   \ar[l]\ar[dr]^-{(f)} \ar[dd] _-{(b)}& \\
       \txt{poly-$T$-complexes} \ar @{<->} [d] _-{(e)}      & &
         \txt{globular\\ $\omega$-groupoids} \ar[dl]^-{(c)} \\
   \txt{simplicial \\ $T$-complexes} \ar[r] _-{(d)}  &\txt{crossed \\ complexes}
   \ar[l]\ar[uu]\ar[ur]
   &
}}\end{equation} Each arrow here denotes an explicit functor which
is an equivalence of categories. The equivalences (a) and (b) are in
\cite{BHig:algcub}; (a) is an essential technical tool in the use of
cubical \ogpd s.  The equivalence (c) is in \cite{BH81:inf}, and
this with (b) implies the equivalence (f); a direct form of this
equivalence is given in the much harder category case in \cite{ABS}.
The equivalence (d) is due to Ashley in \cite{As78}. The equivalence
(e) is due to Jones \cite{Jones}. The different forms of algebra
reflect different geometries, those  of disks, globes, simplices,
cubes, as shown in dimension 2 in the following diagram.
\begin{center}
\begin{overpic}[scale=0.7]{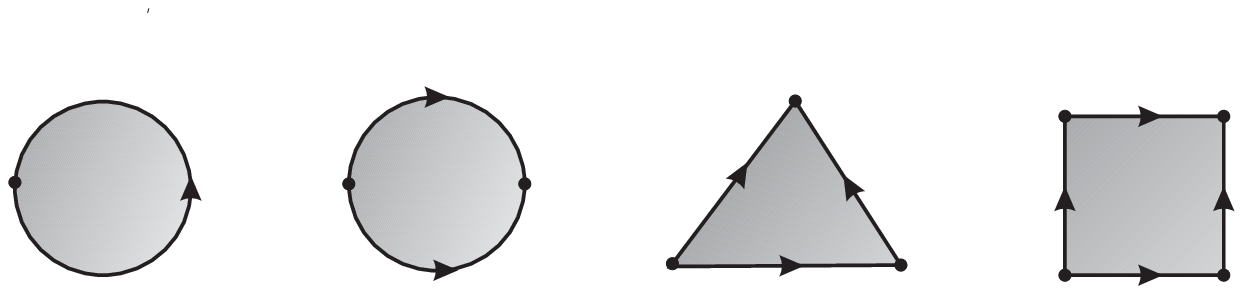}
\end{overpic}
\mbox{ }\\
\mbox{ }\\
\end{center}
It is because the geometry of convex sets is so much more
complicated in dimensions $> 1$ than in dimension $1$ that  new
complications emerge for the theories of higher order group theory
and of higher homotopy groupoids.

A classical homotopical functor on filtered spaces is the {\it
fundamental crossed complex} $\Pi X_*$ of a filtered space, defined
using relative homotopy groups  (in the case $X_0$ is a singleton)
by Blakers, \cite{Bl48}. Major achievements of the papers
\cite{BHig:algcub,BH81:col} were \begin{itemize}
\item to define a homotopical
functor, which here we call $\rho^\Box$, from filtered spaces to
cubical $\omega$-groupoids with connections (and hence also to
cubical $T$-complexes), which in dimension $n$ is the filter
homotopy classes rel vertices of filtered maps $I^n_* \to X_*$ (but
see Remark \ref{rem:J0});
\item to prove that this functor preserved certain colimits;
\item to relate $\rho^\Box$ with the classical functor $\Pi$ from
filtered spaces to crossed complexes, and so to prove that $\Pi$
preserves certain colimits.
\end{itemize}
The proofs {\it do not involve traditional techniques such as
singular homology or simplicial approximation.} The results give
nonabelian information in dimensions $\leq 2$, and in higher
dimensions give information on the action of the fundamental group.
Thus the Relative Hurewicz Theorem is a corollary of a HHvKT
\cite{BH81:col}. Analogous  methods were used by Ashley in
\cite{As78} to define a functor $\rho^\Delta$ from filtered spaces
to simplicial $T$-complexes, and his ideas contributed to
\cite{BH81:col}.

However there has been a lack of a directly defined homotopical
functor from filtered spaces to globular $\omega$-groupoids, and
this gap we will fill in this paper.

The definition of classifying space is most convenient  via well
developed simplicial  constructions. In this way we get the
classifying space of a crossed complex, \cite{BH91}. Its properties
are further exploited in, for example,
\cite{BGPT-equivariantII,martins,martinsporter}.

The equivalence of the category of globular \ogpd s with the
category of cubical \ogpd s with connection, and the monoidal closed
structure on the latter constructed in \cite{BH87}, implies a
monoidal closed structure on the category of globular \ogpd  s.
Further  it is shown in \cite{BH91} that the simple rule $[f]\otimes
[g] \mapsto [f \otimes g]$ gives  a natural transformation
$$ \rho^\Box \, X_* \otimes \rho^\Box \,Y_* \to \rho^\Box (X_* \otimes
Y_*)$$ for any filtered spaces  $X_*, Y_*$, where $X_* \otimes Y_*$
is the usual tensor product of filtered spaces given by $$(X_*
\otimes Y_*)_n = \bigcup _{p+q=n} X_p \times Y_q.$$  The induced
transformation on crossed complexes is shown in \cite{Ba-B} to be an
isomorphism if $X_*, Y_* $ are cofibred and connected. It follows
from the above that there is a natural transformation
$$ \rho^\bigcirc \,X_* \otimes \rho^\bigcirc \,Y_* \to \rho^\bigcirc (X_* \otimes
Y_*).$$ This could be difficult to construct directly. This natural
transformation may be used to enrich the category of filtered spaces
over the monoidal closed category  \ogpd s of globular
$\omega$-groupoids.

It should be apparent from the above that it is the cubical case
which gives the possibility of formulating and of proving theorems;
the basic reason is that cubical theory is handy for subdivision and
its inverse, multiple compositions, and is also good for tensor
products. Many theorems  can then, by equivalences of categories, be
translated to the other cases. However the proofs for the cubical
cases, particularly the properties of  thin elements and
$T$--complexes,  involve also the use of crossed complexes and the
equivalence of categories (a,b) of diagram \eqref{diag:manycats}.
Crossed complexes also have a well developed homotopy theory,
\cite{BG:model}, and have a clear relation with chain complexes with
operators, \cite{BH90}. The relation with simplicial theory is
useful because of the wide development of simplicial theory.
Finally, the relation with the globular theory could be useful
because of the wide familiarity of uses of weak structures and lax
functors and natural transformations: for example, compare the
discussion of Schreier theory using  crossed complexes in
\cite{BH82:ext,BP} with the use of 2-groupoids in
\cite{BBF-schreier}. Calculational applications are usually made
using crossed complexes. For example, the paper \cite{BP}  uses the
notion of small free crossed resolution to give small descriptions
of some nonabelian extensions.

\section{Disks, globes,  and cubes} \label{sec:disksglobes}
Our results follow from an  analysis of the relations between globes
and cubes. These results are probably well known but need to be done
carefully for our purposes.

We give  real  space $\R^n$ the Euclidean norm $\|x\|^2= x_1^2 +
x_2^2 + \cdots+ x_n^2$. We embed $\R^n$ in $\R^{n+1}$ as usual by $x
\mapsto (x,0)$.  The $n$-cube $I^n$ will be the subset of $\R^n$ of
points $x$ such that $|x_i|\leq 1$ for all $i$. Thus $I=I^1$ is
identified with $[-1,1]$ and we also identify $I^n$ with the
$n$-fold product of $I$ with itself.

The {\it $n$-disk} is the subspace $D^n$  of $\R^n$ of points $x$
with $\|x\|\leq 1$. The  $(n-1)$-sphere $S^{n-1}$ is the subspace of
$D^n$ of points $x$ with $\|x\|=1$.

We  define the {\it $n$-globe} $G^n$ to be $D^n$ as a space, but
with the cell structure $$G^n= e^0_\pm \cup e^1_\pm \cup \cdots \cup
e^{n-1}_\pm \cup e^n.$$ Here for $i < n$ the closed cell
$\bar{e}^i_\pm $ is the set of points $x=(x_1, \ldots, x_n) \in G^n$
such that $\|x\|=1,\,  x_j=0$ for $j < n-i$ and $\pm x_{n-i}\geq 0$.
This convention is in keeping with the relationship with cubes which
we find convenient. Note that the $(n-1)$--skeleton of $G^n$ is
contained in $S^{n-1}$.

For each of $Q=\Delta, \Box, \bigcirc$ we have a  singular complex
$S^Q \,X$ of a topological space $X$, giving the well known
simplicial and cubical singular complex, and also a `globular'
singular complex consisting of maps $G^n \to X$. We will later
describe this as a `globular set'.
\begin{Def}
We now define by induction maps $\phi_n: I^n \to G^n, n \geq 1 $,
with the following properties, for $x=(x_1, \ldots,x_n)\in I^n$:
\begin{enumerate}[(i)]
\item $\phi_1(x_1)=x_1$;
\item  $|x_i|=1$ for some $i=1, \ldots, n$ if and only if   $
\|\phi_n(x)\|=1$;
\item $|x_i|=1 $ for some $i=1, \ldots, n$  implies $(\phi_n(x))_j=0 $
for $j <i$.
\end{enumerate}
We set for $x=(t,y) \in \R \times \R^{n-1}$:\begin{equation}
\phi_n(t,y)= ( t\sqrt{1-\|\phi_{n-1}(y)\|^2}, \phi_{n-1}(y)).
\end{equation}First note that if $x=(t,y)$ then
$$\|\phi_n(x)\|^2= t^2 +(1-t^2)\|\phi_{n-1}(y)\|^2.$$
This easily proves (ii) and (iii) by induction. \hfill $\Box$
\end{Def}

 The maps $\phi_n: I^n \to G^n$ induce a map
$\bar{\phi}:S^\bigcirc\,X \to S^\Box\,X$.

We define  the  {\it globular subset } $\gamma K$ of a cubical set
$K$  to agree with $K$ in dimensions $0,1$ and to be in dimension $n
\geq 2$ the set of $k$ such that $\partial^\pm _i k \in \Im
\eps_1^{i-1}, \,i=2, \ldots ,n$.

\begin{prop}
The image of $\bar{\phi}:S^\bigcirc\,X \to S^\Box \,X$ is exactly
the globular subset of $S^\Box\,X$.
\end{prop}
\begin{proof}We prove by induction from the formula for
$\phi_n$ that the image is globular. Let $p_1^i: \R^{n}\to \R^{n-i}$
be the projection omitting the first $i$ coordinates.   Suppose that
$\phi_{n-1}\bar{\partial}^\pm_i=f_{n-1}p_1^{i-1}$. Then
$\phi_{n}\bar{\partial}^\pm_{i+1}=f'_{n-1}p_1^{i}$ where $
f'_{n-1}(x)= (0,f_{n-1}(x))$.

For the converse, we prove by induction that these are the only
identifications that $\phi_n$ makes. Suppose
$\phi_n(t,y)=\phi_n(t',y')$. Then $\phi_{n-1}(y)= \phi_{n-1}(y')$
and $$t\sqrt{1-\|\phi_{n-1}(y)\|^2} =
t'\sqrt{1-\|\phi_{n-1}(y')\|^2}.$$ Thus if $\|\phi_{n-1}(y)\|\ne 1$
then $t=t'$. But $\|\phi_{n-1}(y)\| = 1$ implies some $|y_i|=1$, by
the inductive hypothesis.
\end{proof}

Let $X_*$ be a filtered space. Then we obtain three {\it filtered
singular complexes} $R^Q\,X_*$ for $Q=D,\bigcirc,\I$ defined as
graded sets by
$$(R^Q\,X_*)_n= \FTop(Q^n_*,X_*). $$ There are  also associated graded homotopy
sets $ \rho^Q\,X_*$ which in dimension $n$ are given by  the
quotient maps
$$p^Q:R^Q\,X_*\to\rho^{\! Q} \,
X_*= R^Q \, X_*/\equiv$$ where $\equiv$ is the relation of homotopy
rel vertices through filtered maps.

In the cases $Q=D,\Box$ it is known that these graded sets obtain
additional structure giving us for $Q=D$ the fundamental crossed
complex $\Pi\,X_*$ and for $Q=\Box$ what is called in
\cite{BH81:col} the {\em fundamental (cubical) $\omega$-groupoid
(with connections)} of $X_*$. However the proof that the standard
compositions on $R^\Box \,X_*$ are inherited by $\rho^\Box \,X_*$ is
non trivial, as is the crucial result that $p^\Box$ is a Kan
fibration of cubical sets.
\begin{rem}\label{rem:J0}
  In \cite{BH81:col}, the homotopies are not taken rel vertices and
  a condition $J_0$ is imposed, that each map $\dot{I}^2 \to X_0$  may
  be extended to a map $I^2 \to X_1$. This condition is in many ways
  inconvenient. The filling processes used in the proofs can all
  be started by assuming instead that the homotopies are rel vertices so that the
  maps $\dot{I}^2 \to X_0$ required to be extended are in fact all
  constant. The details will be available in \cite{BHSbook}. \hfill $\Box$
\end{rem}

Our first main result is:
\begin{thm} The induced map $$\phi^*:\rho^\bigcirc  \, X_* \to
\rho^\Box \,X_*$$ is injective.
\end{thm}
\begin{proof}Let $[\alpha],[\beta] \in (\rho^\bigcirc\, X_*)_n$ be such that
$\phi^*[\alpha]=\phi^*[\beta]$, that is $[\alpha\phi]=[\beta\phi]$
in $(\rho^\Box \,X_*)_n$. Let $H:\alpha\phi \equiv\beta\phi$ be such
a homotopy. Then $H$ is a map $I^{n+1} \to X$ such that writing
$I^{n+1}=I^n \times I$, each $H_t: I^n \to X$ is a filtered map.

We  use a folding map $\Phi : I^n \to I^n$ given  by Definition 3.1
of \cite{ABS} (see Definition \ref{def:fold}) which has the property
that $\Phi$ factors through $\phi$.

We now define a new homotopy $K_t=\Phi H_t: I^n \to X$.  Then $K_t$
is a globular homotopy $\Phi \alpha \phi \equiv \Phi \beta \phi $.
But, by assumption, $\alpha \phi, \beta \phi $ are already globular
maps. So the proof is completed with the following lemma.
\begin{lem}
If $a: I^n_* \to X_*$ is a globular map, then $\Phi a$ is globularly
equivalent to $a$.
\end{lem}
\begin{proof}
Since $\Phi$ is a composition of the folding operations $\psi_i$, it
is sufficient to prove that $\psi_i a \equiv_\bigcirc a$. We follow
the proof of \cite[Proposition 3.4]{ABS}. By  the definition of
$\psi_i$:
$$  \psi_i a =
\Gamma^+_i \partial^-_{i+1} a \circ_{i+1} a \circ_{i+1}\Gamma^-_i
\partial^+_{i+1} a.$$ But $\partial^-_{i+1} a$ and
$\partial^+_{i+1} a$  By the laws (A2) we obtain, since $a$ is
globular, that
$$\Gamma^{-\pm}_i \partial^\pm _{i+1}a\in \Im \Gamma ^{-\pm}_i
\epsilon _i = \Im \epsilon_i^2= \Im \epsilon _{1+1}\epsilon _i.$$ So
standard contractions of the two cubes $\Gamma^{-\pm}_i \partial^\pm
_{i+1}a$ yield a  homotopy of $\psi_i a \equiv_\bigcirc a$ through
globular maps.
\end{proof}
It now follows that $\alpha , \beta : G^n_* \to X_* $ are globularly
equivalent.
\end{proof}
 This proof is a higher
dimensional version of an argument in section 6 of
\cite{BHKP:double}.

\begin{cor} The compositions in $\rho^\I \,X_*$ are inherited by
$\rho^\bigcirc\,X_*$ to give the latter the structure of globular
$\omega$-groupoid.
\end{cor}
We do not know how to prove directly that $\rho^\bigcirc\, X_*$ may
be given this structure of globular $\omega$--groupoid.

\section{The free globular $\omega$-groupoid on one generator}
Let $X_*$ be a filtered space. Then we have a diagram of maps of
homotopy sets
\begin{equation}
  (\Pi X_*)_n \labto{i} (\rho^\bigcirc X_*)_n \labto{j}
  (\rho^ \Box X_*)_n.
\end{equation}
We know from \cite{BH81:col} that the composition $j\circ i$ is
injective. We already know that $j$ is injective. It follows that
$i$ is injective. Thus the globular $\omega$-groupoid
$\rho^\bigcirc X_* $ contains the crossed complex $\Pi X_*$, and the
results of \cite{BH81:inf} show that the latter generates the former
as \ogpd .

We need below the following result.
\begin{thm} \label{realise} If $G$ is a globular $\omega$-groupoid, then there is a
filtered space $X_*$ such that $\rho^\bigcirc\,X_* \cong G$.
\end{thm}
\begin{proof}Let $C$ be the crossed complex associated with the $\omega$-groupoid $G$
 under the equivalence (c) of diagram \eqref{diag:manycats}. By
Corollary 9.3 of \cite{BH81:col}, there is a filtered space $X_*$
such that $\Pi \,X_* \cong C$. (Here $X$ is the classifying space
$BC$ filtered by $X_n= BC^{(n)}$ where $C^{(n)}$ is the $n$th
truncation of $C$.)  It follows that $\rho^\bigcirc \,X_* \cong G$.
\end{proof}

\begin{thm}
The globular $\omega$--groupoid  $\rho^\bigcirc G^n_* $ is the free
globular $\omega$--groupoid on the class of the identity map, and
its associated crossed complex is isomorphic to $\Pi G^n_*$.
\end{thm}
\begin{proof}
Let $\iota: G^n_* \to G^n_*$ denote the identity map, and $[\iota]$
its class in $\rho^\bigcirc G^n_* $. Let $H$ be a globular
$\omega$-groupoid and let $x \in H_n$. We have to show there is a
unique morphism $\alpha:\rho^\bigcirc G^n_* \to H$ such that
$\alpha[\iota]=x$. By Theorem \ref{realise} we may assume $H$ is of
the form $\rho^\bigcirc X_* $ for some filtered space $X_*$. Then
$x$ has a representative $g:G^n_* \to X_*$. It follows that
$\rho^\bigcirc(g)([\iota])= x$. This proves existence of such a
morphism.

Suppose $\beta:\rho^\bigcirc G^n_* \to H$ is another morphism such
that $\beta([\iota])=x$. Then $\gamma(\alpha),\gamma(\beta): \Pi
G^n_*$ agree on the generating element $c^n \in
\pi_n(G^n,G^n_{n-1},1)$ of that group. However $\Pi G^n_*$ is
generated as crossed complex by all elements $\Phi dc^n \in
\pi_r(G^n_r,G^n_{r-1},1)$ for all globular face operators $d$ from
dimension $n$ to dimension $r$ for $0 \leq r \leq n$. Since $\alpha,
\beta$ are morphisms of \ogpd s, $\alpha(\Phi dc^n) = d\alpha\Phi
c^n= d\beta\Phi c^n=\beta (d\Phi c^n)$ . Therefore $\alpha, \beta$
agree on $\Pi G^n_*$. But the latter generates $\rho^\bigcirc G^n_*$
as \ogpd . So $\alpha=\beta$.
\end{proof}
The form of this crossed complex may be deduced from the cubical
Homotopy Addition Lemma, \cite[Lemma 7.1]{BHig:algcub}.
\[ \delta x = \begin{cases}- x^+ _1 - x^-_2 + x^-_1 + x^+ _2 &
 \text{if }   n = 2, \\ -  x^+ _3 - ( x^-_2)^{u_{2}\x} -
x^+ _1 + ( x^-_3)^{u_{3}\x} +  x^+ _2
+ ( x^-_1)^{u_{1}\x} & \text{if } n = 3, \\
 \sum_{i=1}^{n} (-1)^i \{  x^+ _i - ( x^-_i)^{u_{i}\x} \} &
 \text{if } n \geq 4
\end{cases}
 \](where $u_i = \partial^+ _1 \partial^+ _2 \cdots
\widehat{\imath} \cdots
 \partial^+ _{n+1})$. In the case when $x$ is globular, this reduces
 to
 $$\delta x = - x^+ _1+ x^-_1  \text{ if } n\geq 2.
$$ Notice that this is a groupoid formula if $n=2$.

\section{Closed monoidal structure}
The category  of cubical \ogpd s with connection is monoidal closed,
\cite{BH87}.  We recall from that paper how the tensor product is
defined.

For cubical \ogpd s $F,G,H$, we  define a {\it bimorphism}
\begin{equation}
  b: F,G \to H
\end{equation}
 to be a family of functions $b= b_{p,q}:
F_p \times G_q \to H_{p+q}$ such that if $ x \in F_p, \; y \in G_q $
and $ p+q =n$ then:
\begin{enumerate}[(i)]
\item $\del^\alpha _ib(x , y)=
\begin{cases}b(\del^\alpha _i x , y) & \text{if } 1 \le i \le p,\\
              b(x , \del^\alpha _{i-p}y)& \text{if } p+ 1 \le i \le n;
\end{cases}$
\item $\eps_ib(x , y) =
\begin{cases}b(\eps_i x , y) & \text{if } 1 \le i \le p+1,\\
              b(x , \eps_{i-p}y)& \text{if } p+ 1 \le i \le n+1;
\end{cases}$
\item $\Gamma_ib(x , y) =
\begin{cases}b(\Gamma_i x , y) & \text{if } 1 \le i \le p,\\
              b(x , \Gamma_{i-p}y)& \text{if } p+ 1 \le i \le n;
\end{cases}$
\item $b(x\circ_i x' ,y) =
b(x , y)\circ_i b(x' , y) \text{ if } 1 \le i \le p \text{ and }
x\circ_i x'\text{ is defined in } F$;
\item $b(x , y\circ_j y') =
b(x , y)\circ_{p+j} b(x , y') \text{ if }  1 \le j \le q \text{ and
} y \circ_j y' \text{ is defined in } G$;
\end{enumerate}

The {\it tensor product} of cubical $\omega$--groupoids  $F,G$  is
given by the the universal  bimorphism $F,G \to F \otimes G$: that
is any bimorphism $F,G \to H$ uniquely factors through a morphism $F
\otimes G \to H$.

We next recall  a result from \cite{BH91}.
\begin{prop}
Let $X_*,\,Y_*$ be filtered spaces. Then there is a natural
transformation $$ \eta: \rho^\Box\, X_* \otimes \rho^\Box
\,Y_*\to\rho^\Box (X_* \otimes Y_*).$$
\end{prop}
\begin{proof}
This natural transformation is determined  by the bimorphism
   \begin{equation*}
    ([f],[g]) \mapsto [f \otimes g]
  \end{equation*}
where $f: I^p_* \to X_*, g: I^q_* \to Y_*$. The proof that this is
well defined and  gives a bimorphism is routine, given the geometry
of the cubes, that $I^p_* \otimes I^q_* \cong I^{p+q}_*$, and the
well definedness of compositions on filter homotopy classes, as
proved in \cite{BH81:col}.
\end{proof}
It is proved in \cite{BH91}, by considering the corresponding free
crossed complexes, that this morphism is an isomorphism if $X_*,Y_*$
are skeletal filtrations of $CW$-complexes, and in \cite{Ba-B} that
this is an isomorphism if $X_*,\,Y_*$ are connected and cofibred.

Because the categories of cubical and of globular  are equivalent,
and the former has a monoidal closed structure, this is inherited by
the latter.

So we deduce from the above results:
\begin{thm} Let $X_*,\,Y_*$ be filtered spaces. Then there is a
natural transformation $$ \eta: \rho^\bigcirc\,X_* \otimes
\rho^\bigcirc\,Y_*\to\rho^\bigcirc(X_* \otimes Y_*)$$ which is an
isomorphism if $X_*,\,Y_*$ are connected and cofibred.
\end{thm}

\section{The higher homotopy van Kampen Theorem}
Suppose for the rest of this section that $X_*$ is a filtered space.
We suppose given a cover ${\cal{U}} = \{ U^\lambda \}_{\lambda \in
\Lambda}$ of $X$ such that the interiors of the sets of ${\cal{U}}$
cover $X.$ For each $\zeta \in \Lambda^n$ we set $U^\zeta =
U^{\zeta_{1}} \cap \cdots \cap U^{\zeta_{n}}, U^\zeta_i = U^\zeta
\cap X_i.$  Then $U^\zeta_0 \subset U^\zeta_1 \subset \cdots$ is
called the {\it induced filtration} $U^\zeta_*$ of $U^\zeta$. So the
globular homotopy $\omega$-groupoids in the following
$\varrho^\bigcirc$-{\it diagram} of the cover are well defined:
\begin{equation} \xymatrix{\bigsqcup_{\zeta \in \Lambda^{2}}
\varrho^\bigcirc \,U^\zeta_* \ar @<0.5ex>[r]^a \ar @<-0.5ex>[r]_b &
\bigsqcup_{\lambda \in \Lambda} \varrho^\bigcirc\, U^\lambda_* \ar
[r] ^-{c} & \varrho^\bigcirc \,X_*}
\end{equation}Here $\bigsqcup$ denotes disjoint union (which is the same as
coproduct in the category of globular $\omega$-groupoids); $a, b$
are determined by the inclusions $a_\zeta : U^\lambda \cap U^\mu
\rightarrow U^\lambda , b_\zeta : U^\lambda \cap U^\mu \rightarrow
U^\mu$ for each $\zeta = (\lambda , \mu ) \in \Lambda^2$; and $c$ is
determined by the inclusions $c_\lambda : U^\lambda \rightarrow X.$

\begin{Def}\label{def:connfilt}
A filtered space $X_*$ is said to be {\it connected}
if the following conditions hold for each $n \geqslant 0:$ \\
$\bullet$ If $r > 0,$ the map $\pi_0 X_0 \rightarrow \pi_0 X_r,$
induced by inclusion, is surjective;  i.e. $X_0$ meets
all path connected components of all stages of the filtration $X_r$. \\
$\bullet$ (for $n \geqslant 1$): If $r > n$ and $x \in X_0$, then $
\pi_n (X_r , X_n , x) = 0.$\hfill $\Box$ \end{Def}

\begin{thm} \label{thm:hvkt}  Suppose that for
every finite intersection $U^\zeta$ of elements of ${\cal{U}}$, the
induced filtration $U^\zeta_*$ is connected. Then
\begin{enumerate}
  \item[\rm{(C)}] $  X_*$ is connected;
  \item[\rm{(I)}] $c$ in the above $\varrho^\bigcirc$--diagram   is the coequaliser
           of $a, b$ in the category  of globular $\omega$--groupoids.
\end{enumerate}
\end{thm}
\begin{proof}
This follows from Theorem B of \cite{BH81:col}, i.e. the analogous
theorem for $\rho^\Box$, and the fact that the equivalence from the
category of globular $\omega$-groupoids to that of cubical
$\omega$-groupoids with connections takes $\rho^\bigcirc\,X_*$ to
$\rho^\Box \,X_*$.
\end{proof}

\section{Nerves and classifying spaces of globular
$\omega$-groupoids} Here we just show how to define a simplicial
nerve $N^\Delta G$ of a globular $\omega$-groupoid  $G$,  by the
standard procedure:
\begin{equation}
(N^\Delta G)_n= \omega\text{-}\mathsf{ Gpd}(\rho^\bigcirc
\,\Delta^n_*, G) .
\end{equation}
The geometric realisation of this simplicial set then defines the
{\it classifying space} $BG$ of $G$. However it is not so easy to
see how to exploit this. The classifying space of a crossed complex
is applied in for example
\cite{BH91,BGPT-equivariantII,martins,martinsporter}.

\appendix

\section{The globular site}\label{sec:globular}
We now recall  from \cite{BH81:inf} a definition which in
\cite{Street:simplices}, and later work,  is termed that of  a {\it
globular set}. This is a sequence $(S_n)_{n \geq 0}$ of sets with
two families of functions
\begin{alignat*}{2} d^\pm_i:
S_n & \to S_i, &&i=0, \ldots, n-1, \\
s_i: S_i & \to S_n,\; &&i=0, \ldots, n-1,
\end{alignat*}
satisfying the following laws, where $\alpha,\beta=\pm$:
 \begin{enumerate}[(i)]
\item $ d^\alpha_i d^\beta_j= d^\alpha _i \text{ for } i < j,
\alpha,\beta= \pm$;
\item $s_j s_i= s_i \text{ for } i < j$;
\item $d^\beta_j s_i = \begin{cases}
  s^\beta_j & \text{ for } j < i, \\
  1 & \text{ for } j=i, \\
  s_i & \text{ for } j > i.
\end{cases}$
\end{enumerate}

A {\it globular site} $GS$ is a small category such that globular
sets can be identified with contravariant functors $GS \to \Set$. We
want to identify such a site whose objects are the globes $G^n$ of
section \ref{sec:disksglobes}. We therefore define maps
\begin{alignat}{2}
\bar{d}^\pm_i: G^i &\to G^n,  &\bar{s}_i:G^n&\to G^i\\
x &\mapsto (0_{n-i}, \pm \sqrt{1-\|x\|^2}, x),\quad & (x_1,
\ldots,x_n) &\mapsto (x_1, \ldots,x_i)
\end{alignat}
for $i <n$, where $0_j=\underset{j}{\underbrace{(0, \ldots,0)}}$.

\section{The cubical site}\label{sec:cubical}
Let  $K$  be a cubical set, that is, a family of  sets  $\{ K_n;n
\ge 0\} $ with face maps  $ \del _i^\A:K_n  \to   K_{n-1} \; (i =
1,2,\ldots,n;\, \alpha = +,-)$  and degeneracy maps
$\epsilon_i:K_{n-1}     \to   K_n \; (i = 1,2,\ldots,n)$ satisfying
the usual cubical relations:

\begin{alignat*}{2}
         \del_i^\A \del_j^\B &=  \del_{j-1}^\B \del_i^{\alpha}
        &&\hspace{-5em}(i<j),
        \tag*{(B.1)(i)} \\
        \epsilon_i \epsilon _j &=  \epsilon _{j+1} \epsilon _i && \hspace{-5em}(i \le j),
          \tag*{(B.1)(ii)} \\
\del ^{\alpha}_i \epsilon _j &=
                  {\begin{cases} \eps_{j-1} \del _i^\A &\hspace{8em} (i<j)  \\
                          \eps_{j} \del _{i-1}^\A & \hspace{8em}(i>j)  \\
                             \mathrm{id} &\hspace{8em} (i=j)
                  \end{cases}} && \tag*{(B.1)(iii)}   \\
\intertext{ We say that  $K$  is a  {\em cubical  set with
connections}  if  it  has additional structure maps (called {\em
connections}) $\Gamma_i^+,\Gamma_i^-  :K_{n-1}     \to    K_n \;
(i = 1,2,\ldots,n-1)$  satisfying the relations:}
       \Gamma_i^\A   \Gamma_j^\B  & =  \Gamma_{j+1}^\B\Gamma_i^\A &&
        \hspace{-5em}  (i  <  j)
       \tag*{(B.2)(i)} \\
         \Gamma_i^\A   \Gamma_i^\A  & =  \Gamma_{i+1}^\A\Gamma_i^\A &&
       \tag*{(B.2)(ii)} \\
       \Gamma_i^\A   \epsilon_j & = {\begin{cases}  \epsilon_{j-1}\Gamma_i^\A &\hspace{8em} (i < j)\\
                                     \epsilon_{j}\Gamma_{i-1}^\A &\hspace{8em}(i > j)
                                  \end{cases}}&& \tag*{(B.2)(iii)} \\
        \Gamma_j^\A \eps_j &= \eps^2_j=\eps_{j+1}\eps_j, &&  \tag*{(B.2)(iv)} \\
        \del^\A _i \Gamma_j^\B &= {\begin{cases} \Gamma_{j-1}^\B\del^\A_i
        & \hspace{8em}(i<j) \\
                              \Gamma_{j}^\B\del^\A_{i-1}  &\hspace{8em} (i> j+1),
                              \end{cases}}&& \tag*{(B.2)(v) }\\
        \del^\A_j\Gamma_j^\A&= \del_{j+1}^\A   \Gamma _j^\A  =  id, && \tag*{(B.2)(vi)} \\
        \del^\A_j \Gamma^{-\A} _j&=
        \del^\A_{j+1}  \Gamma^{-\A} _j  = \epsilon_j  \del^\A_j.
        &&        \tag*{(B.2)(vii)}
\end{alignat*}
The connections are to be thought  of  as  extra  `degeneracies'.
(A degenerate cube of type  $ \epsilon_j  x$  has a pair of
opposite  faces  equal and all other faces degenerate.  A cube of
type  $ \Gamma_i^\A  x$  has a pair  of adjacent faces equal and
all other faces of type $\Gamma_j^\A  y$  or $\epsilon_j y$ .)

The prime example of  a  cubical   set  with  connections  is the
singular cubical complex  $K=S^\Box \,X$  of a space   $X$.   Here
$K_n$ is the set  of singular $n$-cubes in $X$ (i.e. continuous maps
$I^n \to X$). The face maps are induced as usual by maps
$\bar{\partial}^\pm_i:I^{n-1} \to I^n$ and the degeneracies by the
projections $p_i: I^n \to I^{n-1}$.  The connections $ \Gamma_i^\A
:K_{n-1 } \to K_n$ are induced by the maps $\gamma_i^\A : I^n \to
 I^{n-1}$    defined by
   $$   \gamma _i^\A (t_1 ,t_2 ,\ldots,t_n ) =
             (t_1 ,t_2 ,\ldots,t_{i-1},A(t_i ,t_{i+1}),t_{i+2},\ldots,t_n )
             $$
 where $A(s,t)=\max(s,t), \min(s,t)$ as $\A=-,+$ respectively.

The complex  $S^\Box \,X$  has some further relevant structure,
namely the composition of $n$-cubes in the $n$  different
directions. Accordingly, we define a {\it cubical set with
connections and compositions} to be a cubical set  $K$  with
connections in which each $K_n$ has  $n$ partial compositions $\circ
_j\; (j = 1,2,\ldots,n)$ satisfying the following axioms.

If  $a,b \in K_n$, then  $a\circ _j b$  is defined if and only if
$\del^-_j   b = \del^+_j  a$  , and then
\begin{equation} \begin{cases} \del^-_j  (a\circ _j b) = \del^-_ja & \\
                 \del^+_j  (a\circ _j b) = \del^+_jb & \end{cases}
                 \qquad
 \del^\A_i  (a\circ _j b) =  \begin{cases} \del^\A_ja\circ _{j-1}\del^\A_i b &(i<j) \\
                 \del^\A_i  a\circ _j \del^\A_i b& (i>j), \end{cases}
                 \tag*{(B.3)}    \end{equation}

 {\em The  interchange laws}.  If  $i \ne j$  then
\begin{equation}
      (a\circ _i b) \circ _j  (c\circ _i d) = (a\circ _j c) \circ _i  (b\circ _j d)
      \tag*{(B.4)}
\end{equation}
whenever both sides are defined. (The diagram
$$
  \begin{bmatrix}
               a &b \\c&d
  \end{bmatrix}     \quad \directs{i}{j}
$$ will be used to indicate that both sides of the above equation
are  defined and also to denote the unique composite of the four
elements.)

If  $i \ne j$  then
\begin{align*}
           \epsilon_i(a\circ _j b) &= \begin{cases}
           \eps_ia \circ _{j+1} \eps_ib & (i \le j) \\
           \eps_ia \circ _j\eps_ib & (i >j) \end{cases} \tag*{(B.5)} \\
       \Gamma^\A _i (a\circ _j b)& =  \begin{cases}
           \Gamma^\A_ia \circ _{j+1} \Gamma^\A_ib & (i < j) \\
           \Gamma^\A_ia \circ _j\Gamma^\A_ib & (i >j) \end{cases}
           \tag*{(B.6)(i)} \\
       \Gamma^+_j(a\circ _jb)&=  \begin{bmatrix}\Gamma^+_ja & \eps_j a\\
       \eps_{j+1} a & \Gamma^+_j b \end{bmatrix} \quad \directs{j}{j+1}  \tag*{(B.6)(ii)}\\
       \Gamma^-_j(a\circ _jb)&=  \begin{bmatrix}\Gamma^-_ja & \eps_{j+1} b\\
       \eps_{j} b & \Gamma^-_j b \end{bmatrix}\quad \directs{j}{j+1} \tag*{(B.6)(iii)} \\
\end{align*}        These last two equations are the {\it transport
laws}\footnote{Recall from \cite{BS76} that the term {\it
connection} was chosen because of an analogy with path-connections
in differential geometry. In particular, the transport law is a
variation or special case of the transport law for a
path-connection. }.

It is easily verified that the singular cubical complex   $S^\Box
\,X$ of a space $X$   satisfies these axioms if  $\circ _j$   is
defined by
$$
     (a\circ _j b)(t_1 ,t_2 ,\ldots,t_n ) = \begin{cases}
     a(t_1 ,\ldots, t_{j-1},2t_j,t_{j+1} ,\ldots,t_n) &(t_j \le
     \half)\\
      b(t_1 ,\ldots, t_{j-1},2t_j-1,t_{j+1} ,\ldots,t_n) &(t_j \ge
     \half)\\
     \end{cases}
 $$
whenever   $\del^-_j  b =  \del^+_j  a$.

We will now describe two graded subsets of a cubical set $K$. The
{\it globular subset} $K^\bigcirc$ consists in dimension $n$ of the
elements $a$ such that   $\partial^\A _i a \in \Im \eps_1^{i-1},
i=1, \ldots,n $. The {\it diskal subset} $K^D$ consists in dimension
$n$ of the elements $a$ such that $\partial^\A _i a \in \Im
\eps_1^{n-1}$ for $(\A,i) \ne (-,1)$. Clearly $K^D \subs K^\bigcirc
\subs K$.
\begin{prop} If $K$ is a cubical set with compositions, then the
compositions $\circ_i$ are inherited by $K^\bigcirc$ so that if
$d^\A_i : K^\bigcirc _n \to K^\bigcirc _{n-i}$ is defined by $a
\mapsto (\partial^\A_1)^i(a)$, then $K^\bigcirc$ becomes a globular
set with compositions. If further $K$ is a cubical
$\omega$--category (--groupoid), then $K^\bigcirc$ is a globular
$\omega$--category (--groupoid).
\end{prop}

It is proved in \cite{BHig:algcub} that if $K$ is a cubical
$\omega$--groupoid then $K^D$ inherits the structure of crossed
complex, and in \cite{BH81:inf}, see also \cite{ABS},   that
$K^\bigcirc$ inherits the structure of globular $\omega$-groupoid.

A {\it globular $\omega$-category} is a globular set as above with
category structures $\circ_i$ on $S_n$ $ 0 \leq i \leq n-1$ for each
$n \geq 0$ such that $\circ_i$ has $S_i$ as its set of objects and
$D^-_i, D^+_i, E_i$ as its initial, final, and identity maps. These
category structures must be compatible, that is:
\begin{enumerate}[\rm (i)]
\item if $i > j $ and $\alpha=\pm$ then$$D^\alpha_i(x \circ_jy )= D^\alpha_i x \circ_j D^\alpha_i
y,$$whenever the left hand side is defined;
\item  $E_i(x \circ_jy )= E_i x \circ_j E_i y$ in $S_n$ whenever the
left hand side is defined;
\item (The interchange law) if $i \ne j$ then
$$(x \circ_jy) \circ_i (z \circ_jw) = (x \circ_iz) \circ_j (y \circ_jw)$$
whenever both sides are defined.
\end{enumerate}
It is standard to write both sides of the interchange law (when
defined) as \begin{equation*} \begin{bmatrix}   x& y \\z& w\,
\end{bmatrix} \quad \directs{j}{i}
\end{equation*}

\begin{Def} \label{def:fold}{  Let $K$ be a cubical set with connections and compositions.
The {\it folding operations\/} are the operations
$$\psi_i,\Psi_r,\Phi_m\colon K_n\to K_n$$ defined for $1\leq i\leq
n-1$, $1\leq r\leq n$ and $0\leq m\leq n$ by
\begin{gather*}
 \psi_i x
 =\Gamma^+_i\d^-_{i+1}x\circ _{i+1}x\circ _{i+1}\Gamma^-_i\d^+_{i+1}x,\\
  \Psi_r=\psi_{r-1}\psi_{r-2}\ldots\psi_1,\\  \Phi_m=\Psi_1\Psi_2\ldots\Psi_m
 =\psi_1(\psi_2\psi_1)\ldots(\psi_{m-1}\ldots\psi_1).\end{gather*}
} \hfill $\Box$\end{Def}

Note in particular that $\Psi_1$, $\Phi_0$ and~$\Phi_1$ are identity
operations.

Here is a picture of $\psi_1: K_2 \to K_2$:
\vspace{1.5in}

\begin{center}
\setlength{\unitlength}{0.1in}
\begin{picture}(0,0)(10,10)
 \put(0,0){\framebox(4,12){$x$}} \put(0,4){\line(1,0){4}}
\put(0,8){\line(1,0){4}}

\put(0,1){\line(1,0){3}}

 \put(0,2){\line(3,0){2}} \put(2,4){\line(0,-2){2}}
\put(0,3){\line(3,0){1}} \put(2,4){\line(0,-2){2}}
\put(1,4){\line(0,-2){1}} \put(3,4){\line(0,-2){3 }}

 \put(1,11){\line(2,0){3}}
\put(1,11){\line(0,-1){3}} \put(1,11){\line(0,-2){2}}
\put(2,10){\line(3,0){2}} \put(2,10){\line(0,-2){2}}
\put(3,9){\line(0,-2){1}} \put(3,9){\line(1,0){1}}
 \put(6,6){\directs{1}{2}}
\put(-7,6){\makebox{$\psi_1(x)= $}}

\end{picture}
\end{center}

\vspace{1in}
\begin{prop} \label{prop:foldtoglob} Let $K$ be a cubical set with connections and compositions. The
`folding'operator $\Phi_n\colon K_n\to K_n$ satisfies $\partial ^\pm
_i \Phi_n x \in \Im\epsilon_1^{i-1}$  for $1\leq i\leq n$ and $x \in
K_n$. That is, $\Im \Phi $ is contained in the globular subset of
$K$.
\end{prop}
This is part of Proposition 3.3(iii) of \cite{ABS}. Note that the
compositions are needed to define $\Phi_n$ but this property of
$\Phi_n$ does not require any axioms on the compositions, but only
the properties (B1), (B2) giving the relations between cubical
operations and connections.

\end{document}